\documentclass[runningheads]{llncs}
\usepackage[T1]{fontenc}
\usepackage{graphicx} 
\usepackage{amsmath}
\usepackage{amssymb}
\usepackage[nocomma]{optidef}
\usepackage{booktabs}
\usepackage{longtable}
\usepackage{multirow}
\usepackage{lscape}
\usepackage[title]{appendix}

\newcommand{\sigpi}{\sigma_{ii^{\prime}}\pi_{jj^{\prime}}}
\newcommand{\p}{^{\prime}}
\newcommand{\pp}{^{\prime \prime}}

\title{QUBO Formulations for MIP Symmetry Detection}
\author{Alexander While and Chen Chen}
\institute{The Ohio State University, Columbus OH, USA \\
\email{\{while.1, chen.8018\}@osu.edu}}  
\authorrunning{A. While and C. Chen}

\begin{document}

\maketitle

\begin{abstract}
Formulation symmetry in mixed-integer programming (MIP) can hinder solver performance by inducing redundant search, but detecting such symmetries is also a significant computational challenge. This paper explores the potential for quantum computing to handle symmetry detection. Quantum is a promising alternative to classical compute, but this emerging technology has limited hardware capacity in terms of input problem size. This paper explores the use of Quadratic Unconstrained Binary Optimization (QUBO) models for symmetry detection, as QUBO is the canonical format for quantum optimization platforms. To help address the input size bottleneck, we develop full, reduced, and decomposed QUBO as well as QUBO-Plus formulations for MIP symmetry detection. Computational experiments on the MIPLIB 2017 benchmark are used to estimate the quantum computing resources needed for practical problems. 
\end{abstract}


\section{Introduction}

Consider a generic Mixed-Integer Programming (MIP) problem of the form
\begin{maxi}|s|{x}{c^Tx}{}{\label{eq:generic-ip}}
    \addConstraint{Ax}{\leq b}
    \addConstraint{x}{\in \mathbb{Z}^p \times \mathbb{R}^{n-p}}
\end{maxi}
where $A$ is an $n \times m$ matrix.  This is a general-purpose algebraic modeling paradigm that enables access to powerful solver software; for example, enabling billions of dollars in savings in power systems operations \cite{carlson2012miso}. 

As MIP is NP-hard in general, there is persistent demand for computational improvements of solvers.  In this paper we explore the potential of quantum computing to accelerate MIP solvers, namely by considering the problem of symmetry detection. Symmetry can slow down solvers, as identical solutions may be repeatedly enumerated \cite{margot2009symmetry,bodi2013algorithms,ostrowski2015modified,pfetsch2019computational}. To mitigate this redundant computation,  one must first identify symmetries within the problem. In this paper, we present Quadratic Unconstrained Binary Optimization (QUBO) formulations for detecting symmetries within an MIP that can be accepted and heuristically solved for on a quantum annealer or other quantum-inspired solvers (see, e.g. \cite{Abbott2019anneal,glover2022quantum1,glover2022quantum2}). 

Our work adds to the literature on hybrid schemes for MIP, which has considered quantum subroutines involving: primal heuristics \cite{svensson2023hybrid}; Benders cut generation \cite{zhao2022hybrid,paterakis2023hybrid}; combinatorial/pure integer subproblems \cite{ellinas2024hybrid,wei2024hybrid,ajagekar2022hybrid,brown2024copositive,woerner2021solving}; and branching \cite{montanaro2020quantum}. This paper is unique in considering the problem of formulation symmetry detection, which is notably suspected to be NP-intermediate (see Section~\ref{sec:symmetry}), i.e. a problem in NP that is neither NP-Hard nor in P. In contrast, the aforementioned hybrid schemes all involve NP-hard problems for quantum computing, which may be beyond the realm of provable quantum advantage (see Section~\ref{sec:quantum}). We note also that efforts to accelerate linear programming with quantum computing (see, e.g. \cite{mohammadisiahroudi2025improvements,nannicini2024fast}) can also be incorporated into MIP solvers in hybrid fashion.

\subsection{Quantum Computing}
\label{sec:quantum}
Quantum computing has substantial potential; notably, Integer Factorization (IF) was proven to be in bounded-error quantum polynomial time (BQP) via Shor's Algorithm \cite{shor1994shors}. IF is conjectured to be NP-Intermediate, hence challenging for classical compute; indeed, the effectiveness of the widely-used RSA encryption relies on the intractability of IF. Quantum computing has clear theoretical limits, however, as BQP is a subset of QIP=PSPACE \cite{jain2011qip}; moreover, it remains an open question whether BQP contains NP or vice versa, though there are reasons to suspect they are incomparable (e.g. \cite{aaronson2010bqp}).  

At present, practice falls rather short of theory due to the many challenges with scaling quantum computers to sufficient size. For example, gate-based quantum computers struggle to run Shor's algorithm on numbers with more than 2 digits \cite{willsch2023large}. In contrast, classical methods can handle 240-digit numbers \cite{boudot2020comparing}. Indeed, at present, Shor's algorithm itself seems more effectively realized on a classical computer: Willsch et al. \cite{willsch2023large} recently simulated Shor's Quantum Factoring Algorithm with classical computing to factor a 12-digit number.

    An alternative to gate-based quantum computers are quantum annealers, which are composed of physical qubits arranged in a weighted planar graph. Annealers attempt to determine the values of the qubits that lead to the lowest total energy in the graph. While attaining a true minimum is not guaranteed, the process can be repeated thousands of times very quickly, with each sample taking as little as 20ns \cite{Abbott2019anneal}. This process is easily extended to QUBO formulations, which is the canonical format for quantum-based optimization---a generic setup is provided in Figure \ref{fig:generic-qubo}. This type of quantum computer can thus serve as a fast heuristic in hybrid optimization schemes. We note that quantum annealers cannot run Shor's algorithm, but 7-digit numbers have been factorized via QUBO formulation \cite{peng2019factoring}---up to 13-digits when substantial classical compute is applied for preprocessing \cite{jun2023hubo}.  

There have also been advancements in hybrid quantum-classical hardware and algorithms that enable sampling of solutions to QUBO problems with some constraints: QUBO-Plus models. Common constraints include packing, covering, and knapsack constraints (see \cite{glover2022quantum1,glover2022quantum2} for details). In addition to the aforementioned IF results, there has been substantial recent efforts at benchmarking on optimization problems (e.g. \cite{abbas2024challenges,lubinski2024optimization,proctor2025benchmarking,nannicini2019performance,quinton2025quantum,junger2021quantum}), which tell a similar story: quantum computers may serve as an intriguing, rapid primal heuristic, but they are quite limited in input size at present.

 
\begin{figure}
    \centering
    $$
    \min_{x\in \{0,1\}^n} \sum_{i=1}^n\sum_{j=1}^n Q_{ij}x_ix_j
    $$
    \caption{A generic QUBO formulation}
    \label{fig:generic-qubo}
\end{figure}

\subsection{MIP Symmetry}
\label{sec:symmetry}
The symmetry group of a MIP problem instance is the set of constraint and variable permutations that preserve the feasible region without changing the value of the objective function. However, working with this group is computationally challenging, as even determining feasibility of a MIP instance is NP-hard. Instead, solvers typically focus on \textit{formulation symmetry}: a variable permutation $\pi \in S_n$ such that there exists a corresponding constraint permutation $\sigma \in S_m$ where:
\begin{itemize}
    \item $\pi (\{1,\ldots,p\}) = \{1,\ldots, p\}$ (integer variables preserved),
    \item $\pi(c) = c$ (objective preserved),
    \item $\sigma(b) = b$ (constraint constants preserved), and
    \item $A_{\sigma(i),\pi(j)} = A_{ij}$ (rows permuted by $\sigma$, columns by $\pi$).
\end{itemize}
Formulation symmetry detection can be converted to a graph theoretic problem by representing the MIP formulation as a bipartite vertex- and edge-colored graph. In one partition, vertices represent MIP variables, with colors corresponding to their objective coefficients. In the other partition, vertices represent MIP constraints, with colors representing their constant terms. Edges are colored based on variable coefficients from the constraints. The automorphism group of this graph is equivalent to the formulation symmetry group, and so finding MIP symmetries reduces to finding graph automorphisms; consequently, software such as Nauty \cite{mckay2007nauty} or Bliss \cite{junttila2007engineering} can be applied to determine generators of the group. Once generators are obtained, orbits---that is, sets of variables which can be permuted without changing the MIP formulation---can be calculated and then exploited in the MIP solver. While these graph-based algorithms are often fast in practice, they have worst-case exponential run-times, and so in practice one may terminate search early after finding only a strict subset of generators. Indeed, the problem of graph automorphism detection is suspected to be NP-Intermediate \cite{babai2018group}; however, unlike IF it is a longstanding open question whether the problem is in BQP \cite{aaronson2016space}. 

After running symmetry detection, techniques can be applied to prevent symmetric solutions from being revisited in the branch-and-bound tree. Well-established methods include Isomorphic Pruning \cite{margot2009symmetry}, Orbital Fixing \cite{ostrowski2015modified}, and symmetry-breaking constraints \cite{liberti2012reformulations}. Other attempts have also been made leveraging abstract algebraic tools \cite{bodi2013algorithms}. For descriptions of these techniques, as well as an experimental comparison of their performances, see \cite{pfetsch2019computational}.

Solvers such as SCIP \cite{bolusani2024scip} run their own native implementations of symmetry detection (working directly with MIP formulation data structures) to reduce overhead. In a similar vein, we develop in this paper QUBO formulations specifically catered to formulation symmetry detection. This avoids having to reduce a MIP to existing QUBOs for graph isomorphism \cite{Lucas:2014,calude2017qubo,hua2020improved,wang2025rbm}. 

\section{QUBO Symmetry Detection Formulations}

\subsection{Setup and Notation}
The following setup is used for all our QUBO formulations. The decision variables are square, binary matrices $\pi \in \{0,1\}^{n\times n}$ (for variable permutation) and $\sigma \in \{0,1\}^{m\times m}$ (for constraint permutation).  The output of these variables are interpreted as follows: 
\begin{itemize}
    \item $x_j$ is mapped to $x_{j\p}$ if and only if $\pi_{jj\p}=1$ 
    \item constraint $i$ is mapped to constraint $i\p$ if and only if $\sigma_{ii\p}=1$.
\end{itemize}

The QUBO formulations contain bilinear terms of the form $\sigpi$. Following the aforementioned interpretations, we can see that $A_{ij}$ is mapped to $A_{i\p j\p}$ in the coefficient matrix if and only if $\sigpi=1$. 

Moreover, we define a \emph{reasonable} permutation as:
\begin{itemize}
    \item a permutation of variables that maps variables with the same domain and objective coefficients, or
    \item a permutation of constraints that maps constraints with the same right-hand constants $b$.
\end{itemize}
Accordingly, let $r(a)$ be the index set of variables (resp. constraints) that can be reasonably permuted with $x_a$ (resp. constraint $a$). Furthermore, we denote the set of all reasonable variable permutations as $\Pi := \{\pi_{jj\p}:j\p \in r(j),j\in \{1,\ldots,n\}\}$ with complement $\Pi^\complement$, and the set of all reasonable constraint permutations as $ \Sigma := \{\sigma_{ii\p}: i\p \in r(i),i\in\{1,\ldots,m\}\}$ with complement $\Sigma^\complement$.  Note that $\Pi$ (resp. $\Sigma$) can be seen as a permuted block-diagonal submatrix of $\pi$ (resp. $\sigma$), as illustrated in Section~\ref{sec:example}.

While this definition of \emph{reasonable} is sufficient to ensure symmetry detection for the generic form of MIP seen in (\ref{eq:generic-ip}), we can add modifications to handle MIP structures:
\begin{itemize}
    \item variables with the same upper and lower bounds
    \item constraints with the same constraint sense
\end{itemize}
We can also sharpen the notion of reasonable to reduce formulation size, namely by eliminating:
\begin{itemize}
        \item variables that are present in the same number of constraints
        \item constraints containing the same number of variables
\end{itemize}
These modifications are detailed in Section \ref{sec:experiments}.






\subsection{Full and Reduced QUBOs}
Following in the spirit of QUBO formulations for graph isomorphism, we seek to construct a QUBO with the property that for any optimal solution (with objective value 0) we can easily extract a corresponding formulation symmetry. First, let us ensure that $\pi$ and $\sigma$ are doubly stochastic and thus represent permutations \cite{fiedler1988doubly}. To do this, we incorporate the following penalty terms in the objective:
\begin{align}
    H_{B,\pi} &= \sum_{j=1}^n \left( \sum_{j\p=1}^n\pi_{jj\p}-1\right)^2 + 
    \sum_{j\p=1}^n \left(\sum_{j=1}^n \pi_{jj\p}-1\right)^2 \\
    H_{B,\sigma} &= \sum_{i=1}^m \left( \sum_{i\p=1}^m \sigma_{ii\p} -1 \right)^2 
        + \sum_{i\p=1}^m \left( \sum_{i=1}^m\sigma_{ii\p}-1\right)^2
\end{align}

Note that further along in this paper we will be working with partial matrices, i.e. removing elements of $\pi$ and $\sigma$ from our formulations; for simplicity, we will preserve the notation for these expressions, and which elements are preserved in sums will be clear from context.

Now, let us enforce that only reasonable permutations are allowed. We accomplish this with the following penalty terms:

\begin{align}
    H_{\pi}&= \sum_{\Pi^\complement}\pi^2_{jj\p} \\
    H_{\sigma}&=\sum_{\Sigma^\complement}\sigma^2_{ii\p}
\end{align}

Finally, we penalize any combination of variable and constraint permutations that map values of the coefficient matrix $A$ to entries with different values. 
\begin{align}
    H_A=\sum_{A_{ij} \neq A_{i\p j\p}} \sigpi.
\end{align} 

Putting these together gives our Full QUBO formulation: 
\begin{align}
    \min_{\pi,\sigma} H_{Full}:=H_{B,\pi}+H_{B,\sigma} + H_{\pi}+H_{\sigma}+H_A 
\end{align}
 
This formulation consists of $q_{full}:=n^2+m^2$ variables. This can be reduced by restricting our formulation only to certain variables as defined by the reasonable permutation sets $\Pi$ and $\Sigma$, subsequently allowing us to drop the $H_{\Pi}$ and $H_{\Sigma}$ terms. Hence, dropping certain variables gives the Reduced QUBO:
\begin{align}
    \min_{\Pi,\Sigma}H_{Reduced} := H_{B,\pi}+H_{B,\sigma} + H_A
\end{align}

Note that the variables excluded from the formulation, namely $\Pi^\complement,\Sigma^\complement$, are fixed to zero.  The number of variables in the Reduced QUBO is $q_{Reduced}: = \nu + \mu$, where $\nu =|\Pi|=\sum_{j=1}^n|r(i)|$ and $\mu =|\Sigma|=\sum_{i=1}^m|r(j)|$. Note that $n\leq \nu \leq n^2$ and $m \leq \mu \leq m^2$. The values of $\nu$ and $\mu$ depend on the structure of the MIP, and we explore this empirically with the MIPLIB 2017 problem set in Section~\ref{sec:experiments}. 

\begin{proposition} 
\label{prop:full}
If $H_{Full}^*=0$ for some $(\pi^*,\sigma^*)$, then $\pi^*$ is a formulation symmetry. Moreover, if $\bar \pi$ is a formulation symmetry, then there exists some corresponding $\bar \sigma$ such that $H_{Full}(\bar \pi, \bar \sigma)=0$. 
\end{proposition}

\begin{proof}
Let $H^*_{Full} = 0$. Since $H_{B,\pi}, H_{B,\sigma}, H_{\pi}, H_{\sigma},$ and $H_{A}$ all consist of 
sums of squared terms, then each individual term must be zero. Since $H_{B,\pi}=0$ and $H_{B,\sigma}=0$, $\pi^*$ 
and $\sigma^*$ are doubly stochastic and thus represent permutations. 
$H_{\pi}=0$ and 
$H_{\sigma}=0$ only if reasonable permutations are mapped, which
ensures that $\pi(c)=c$, $\sigma(b)=b$, and integer variables are preserved. Finally, we must have that $A_{\sigma(i),\pi(j)} = A_{ij}$ because $H_A=0$.

Now, suppose $\bar \pi$ is a formulation symmetry and $\bar \sigma$ is a corresponding constraint permutation. Since both matrices are permutations, their matrix representations are doubly stochastic and thus $H_{B,\bar \pi}=H_{B,\bar \sigma}=0$. From the definition of formulation symmetry, we also have that integer variables are preserved, $\bar \pi(c)=c$ and $\bar \sigma(b)=b$, so $H_{\bar \pi}=H_{\bar \sigma}=0$. Finally, since $A_{\bar \sigma(i),\bar \pi(j)} = A_{ij}$, this gives us $H_{A}=0$, and so $H_{Full}^*=0$.

\end{proof}

\begin{corollary}
\label{prop:reduced}
If $H_{Reduced}^*=0$ for some $(\pi^*,\sigma^*)$, then $\pi^*$ is a formulation symmetry. Moreover, if $\bar \pi$ is a formulation symmetry, then there exists some corresponding $\bar \sigma$ such that $H_{Reduced}(\bar \pi, \bar \sigma)=0$. 
\end{corollary}

\begin{proof}
Suppose that $H_{Reduced}^*=0$ for some $(\pi^*,\sigma^*)$. Then by nonnegativity, each term $H_{B,\pi}$ and $H_{B,\sigma}$ must be zero, and so $\pi^*$ and $\sigma^*$ are doubly stochastic (permutation) matrices. Likewise, $H_A=0$. Furthermore, by construction, we have variables in $\Pi^\complement$ and $\Sigma^\complement$ fixed to zero, and so $H_\pi=H_\sigma=0$.  Therefore, $H_{Full}^*=0$ and so $\pi^*$ is a formulation symmetry by Proposition \ref{prop:full}.

Now suppose that $\pi^*$ is a formulation symmetry and $\sigma^*$ is a corresponding permutation of the constraints. Since they are both permutations, their matrix representations are doubly stochastic, so $H_{B,\pi}^*=H_{B,\sigma}^*=0$. By definition of formulation symmetry, the constraint matrix is preserved, thus $H_A^*=0$.  Therefore, $H^*_{Reduced}= H_{B,\pi}^*+H_{B,\sigma}^*+H_A^*=0$. \qed
\end{proof}

\subsubsection{Maximum Required Qubits for DWave Embeddings}
We can estimate the resource requirements of our QUBO formulations. One of the key challenges to solving QUBO problems on a DWave quantum annealer is finding an embedding of the problem on the topology of the hardware, which involves limited number of qubits that are not fully connected. 
The first DWave quantum annealers arranged their qubits along a Chimera graph architecture, which allows efficient embeddings of complete graphs. Subsequent generations involved Pegasus graphs, which built off of the Chimera structure to allow for a higher degree of connectivity \cite{boothby2016fast}.
The next generation of DWave technology will adopt the Zephyr graph, which are described by a grid parameter $g$ and a tile parameter $t$ and denoted as $Z_{g,t}$ \cite{boothby2021zephyr}. Due to physical manufacturing and design concerns, it is easiest to fix the tile parameter to $t=4$ while increasing $g$ in order to increase size of the computer \cite{boothby2016fast}---we will focus on such graphs throughout this section, simply denoted as $Z_g$. 

Using the results in \cite{boothby2021zephyr} and \cite{boothby2016fast}, we can create an upper bound of the number of qubits required to embed our Symmetry Detecting QUBOs. 
\begin{proposition}
    \label{prop:qubits}
    A QUBO problem with $q$ variables can be embedded in polynomial time on a $Z_g$ topology with at least $\frac{(q+8)^2}{8}+q+8$ fully functional qubits.
\end{proposition}

\begin{proof}
    The graph $Z_g$ contains $32g^2+16g$ qubits. As described in \cite{boothby2021zephyr}, the largest complete graph that can be efficiently embedded on $Z_g$ using the algorithm in \cite{boothby2016fast} is of size $16g-8$, all with chains of length $g$.  Therefore, if we need to embed $q$ variables, we need $g \geq \frac{q+8}{16}$ and thus at least $\frac{(q+8)^2}{8}+q+8$ qubits in our Zephyr graph. Once the complete graph is embedded, we then delete the unneeded edges for our particular problem.
\end{proof}

There are a few caveats with the result from Proposition \ref{prop:qubits}. In practice, not all qubits on a quantum annealer are typically operable, and the arrangement of the inoperable qubits could require a larger typology. On the other hand, since the desired graphs are certainly not complete, it is very likely that smaller embeddings could be found. We explore this dynamic in Section \ref{sec:experiments} for our Reduced QUBO.

\subsubsection{QUBO-plus Variations}
We also consider QUBO-plus variants of our formulations, which allows us to transform some of our penalty terms into constraints. This is compatible with quantum computers such as DWave's hybrid solvers, which can handle larger-sized problems of up to 10000 variables \cite{developers2020d}. In particular, we move the requirement that $\pi$ and $\sigma$ are doubly stochastic into linear constraints for both Full and Reduced formulations. Both formulations require the same amount of variables as their respective QUBO formulations with the addition of $2n+2m$ constraints and can be seen in Figure \ref{fig:qubo-plus}. 




\begin{figure}
\noindent \begin{minipage}{0.49\linewidth}
    \begin{align}
   \min_{\pi, \sigma} &\quad H_A +H_{\pi}+H_{\sigma}\\
    \text{s.t.\quad} 
    \sum_{j=1}^n \pi_{jj\p} &=1, \quad j\p \in \{1,\ldots,n\}\\
    \sum_{j\p=1}^n \pi_{jj\p} &=1, \quad j\in \{1,\ldots,n\}\\
    \sum_{i=1}^m \sigma_{ii\p} &=1, \quad i\p \in \{1,\ldots,m\} \\
    \sum_{i\p=1}^m \sigma_{ii\p} &= 1, \quad i \in \{1,\ldots,m\}
    \end{align}
\end{minipage}
\hfill\vline\hfill
\begin{minipage}{0.49\linewidth}

    \begin{align}
    \min_{\Pi, \Sigma} &\quad H_A \\
    \text{s.t.\quad} 
    \sum_{j \in r(j\p)} \pi_{jj\p} &= 1 \quad j\p \in \{1,\dots, n\} \\
    \sum_{j\p \in r(j)} \pi_{jj\p} &= 1 \quad j \in \{1,\dots, n\} \\
    \sum_{i \in r(i\p)} \sigma_{ii\p} &=1, \quad i\p \in \{1,\ldots,m\} \\
    \sum_{i\p \in r(i)} \sigma_{ii\p} &= 1, \quad i \in \{1,\ldots,m\}
    \end{align}

\end{minipage}

\caption{Symmetry Detecting QUBO-Plus Formulations, Full-sized (left) and Reduced (Right)}
\label{fig:qubo-plus}
\end{figure}

\subsection{Decompositions}
To further QUBO formulation size, we consider decomposing over reasonable permutations given by $r(j)$. Our decomposition over $r(j)$ is similar to the Reduced form, however we must take care to ensure that the variables which cannot be reasonably permuted with $x_j$ are fixed to ensure that the constraint matrix remains unchanged. As such, we will still include the variables $\pi_{j\p j\p}$ for all $j\p \notin   r(j)$ and either fix them as 1 or require via penalty terms or constraints that they must equal 1, thus ensuring the position of $x_{j}$ remains fixed. 

Let $H_{F,j} = \sum_{j\p\notin r(j)} (1-\pi_{j\p j\p})^2.$ We denote the complete set of $\pi$ variables in this decomposition as $\Pi_j := \{\pi_{j\p j\pp}: j\p, j\pp \in r(j)\} \cup \{\pi_{j\p j\p}: j\p \notin r(j)\}$. We then have the following QUBO Decomposition over $\Pi_{j}$:

\begin{align}
\min_{\Pi_j,\Sigma} H_{Decomp} = H_{B,\pi}+H_{B,\sigma}+H_A+H_{F,j} 
\end{align}
This formulation involves $q_{Decomp}=|r(j)|^2+\left(n-|r(j)|\right) + \mu$ variables. As with the Reduced form, the size depends on problem structure, which we explore empirically in Section~\ref{sec:experiments}. Note that this formulation is restricted in the sense of only providing symmetries on the set of variables symmetric to $\pi_j$ rather than across the entire MIP. We note also that, in principle, decomposition could also be done over the set of constraints in $r(i)$, although it is less desirable as the end-goal is generally to identify symmetric variables.

\begin{proposition}
    If $H_{Decomp}^*=0$ for some $\pi^*,\sigma^*$, then $\pi^*$ is a formulation symmetry that fixes the position of all $x_{j\p}$ which cannot be reasonably permuted with $x_j$ . Furthermore, if $\bar \pi$ is a formulation symmetry that fixes the position of all $x_{j\p}$ which cannot be reasonably permuted with $x_j$, then there exists some $\bar \sigma$ such that $H_{Decomp}(\bar \pi, \bar \sigma)=0$. 
\end{proposition}

\begin{proof}
    Suppose $H_{Decomp}^*=0$ for some $\pi^*,\sigma*$. Since each term is nonnegative, we must have that each term is zero-valued at $(\pi^*,\sigma^*)$. As $H_{B,\pi}=0$ and $H_{B,\sigma}=0$, then $\pi^*$ and $\sigma^*$ are double stochastic and thus represent permutations.  Since the formulation only contains variables that represent mappings within reasonable permutations, the integer variables, objective coefficients, and constraint constants are preserved. Since $H_{F,j}=0$, the variables which cannot be reasonably permuted with $\pi_j$ must be fixed and thus the permutations remain valid. Finally, if $H_A=0$, then the coefficient matrix is preserved. 

    Now consider some formulation symmetry $\bar \pi$ that fixes the position of all $x_{j\p}$ which cannot be reasonably permuted with $x_{j}$. By definition, there exists some corresponding $\bar \sigma$ that creates a valid permutation of the constraints. Since these are permutations, their matrix representations are doubly stochastic and so $H_{B,\pi}=H_{B,\sigma}=0$.  If the position of $x_{j\p}$ is fixed, then $\bar x_{j\p j\p}=1$ and thus $H_{F,j}=0$.  Finally, formulation symmetries preserve the coefficient matrix, so $H_A=0$.
    \end{proof}

We also have the QUBO-Plus Decomposition over $\Pi_{j}$:
\begin{align}
    \min_{\Pi_j,\Sigma} & H_A \\
    \text{s.t.\quad} 
    \sum_{j\p \in r(j)}\pi_{j\p j\pp}&=1, \quad j\pp \in r(j)\\
    \sum_{j\pp \in r(j)}\pi_{j\p j\pp} &=1, \quad j\p \in r(j)\\
    \sum_{i\in r(i)} \sigma_{ii\p} &=1, \quad i\p \in \{1,\ldots,m\} \\
    \sum_{i\p \in r(i)} \sigma_{ii\p} &= 1, \quad i \in \{1,\ldots,m\} \\
    \pi_{j\p j\p}&=1, \quad j\p \notin r(j)
\end{align}
This again features $q_{Decomp}$ variables as well as $2|r(j)|+2m+\left(n-|r(j)|\right)$ constraints.

\section{Experiments}
\label{sec:experiments}
\subsection{Example}
\label{sec:example}
We begin with a simple knapsack problem as an example MIP instance:

\begin{maxi}|s|{x\in \mathbb{Z}^{7}}{x_1+x_2+x_3+2x_4+2x_5+2x_6+3x_7}{}{}
    \addConstraint{x_1+x_2+2x_3+x_4+x_5+x_6+x_7}{\leq 100}{}    
\end{maxi}

\begin{figure}
    \centering
    
    \includegraphics[width=0.3\textwidth]{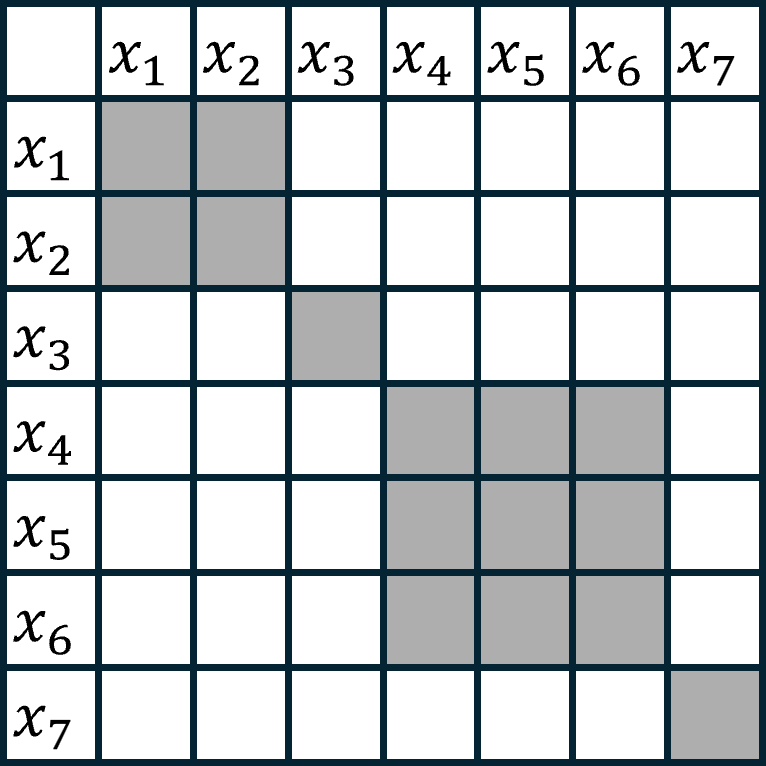}
    \includegraphics[width=0.3\textwidth]{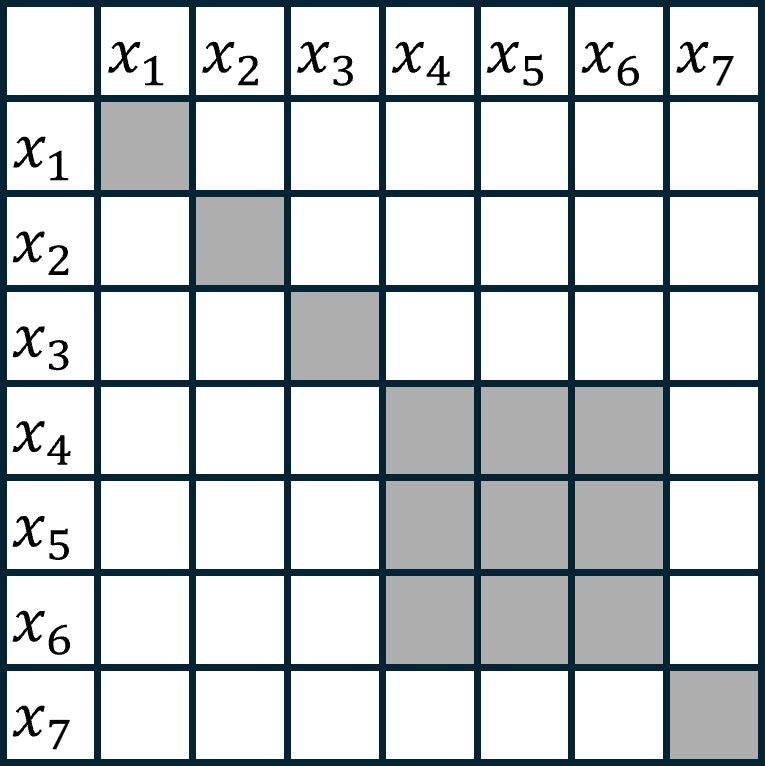}
    \caption{Visualizations of the Reduced (left) and Decomposed (right) forms.}
    \label{fig:visualizations}
\end{figure}

The problem has the following orbits: 
$\{\pi_1,\pi_2\},\{\pi_3\}, \{\pi_4,\pi_5,\pi_6\},\{\pi_7\}$. In our full QUBO, we have decision variables $\pi \in \{0,1\}^{7\times 7}$ and $\sigma_{1,1}$. Thus, $q_{Full} = 7^2+1^1=50$. We can, however, reduce the problem based on the following sets of reasonable permutations:
\begin{align*}
    r(1) &= \{1,2\} \\
    r(2) &= \{1,2\} \\
    r(3) &=\{3\} \\
    r(4) &= \{4,5,6\} \\
    r(5) &= \{4,5,6\} \\
    r(6) &= \{4,5,6\} \\
    r(7) &= \{7\} 
\end{align*}
The $\pi$ variables that are excluded are shown in Figure \ref{fig:visualizations}. 
$\nu =15$ while $\mu=1$, so $q_{Reduced}=15+1=16$. We can quantify the effectiveness of the reduction by calculating $\frac{q_{Full}}{q_{Reduced}} = 0.32$, so the Reduced form requires 32\% as many variables as the Full. Should we require our QUBO to be even smaller, we can use the Decomposed form over any of the sets of reasonable permutations, say $\Pi_4$, at the potential cost of dropping certain symmetries. Again, we can visualize which variables are excluded by looking at Figure \ref{fig:visualizations}. Here we now have $q_{Decomp}=13$, requiring 26\% of the variables of the Full formulation. 



\subsection{Experiments with MIPLIB 2017}
\begin{figure}
    \centering
    \includegraphics[width=0.49\linewidth]{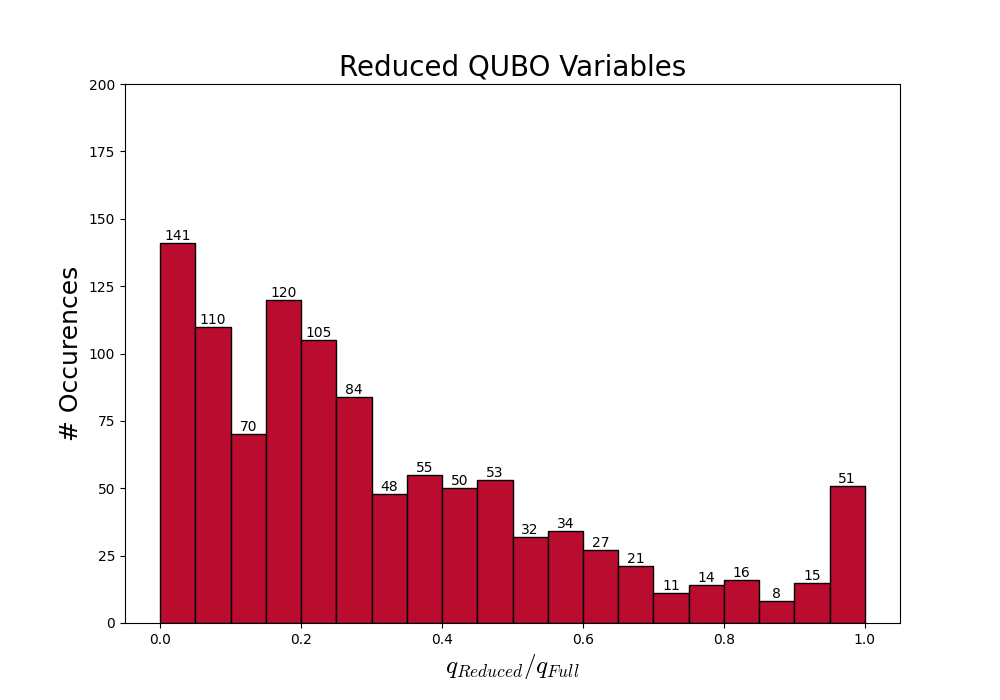}
    \includegraphics[width=0.49\linewidth]{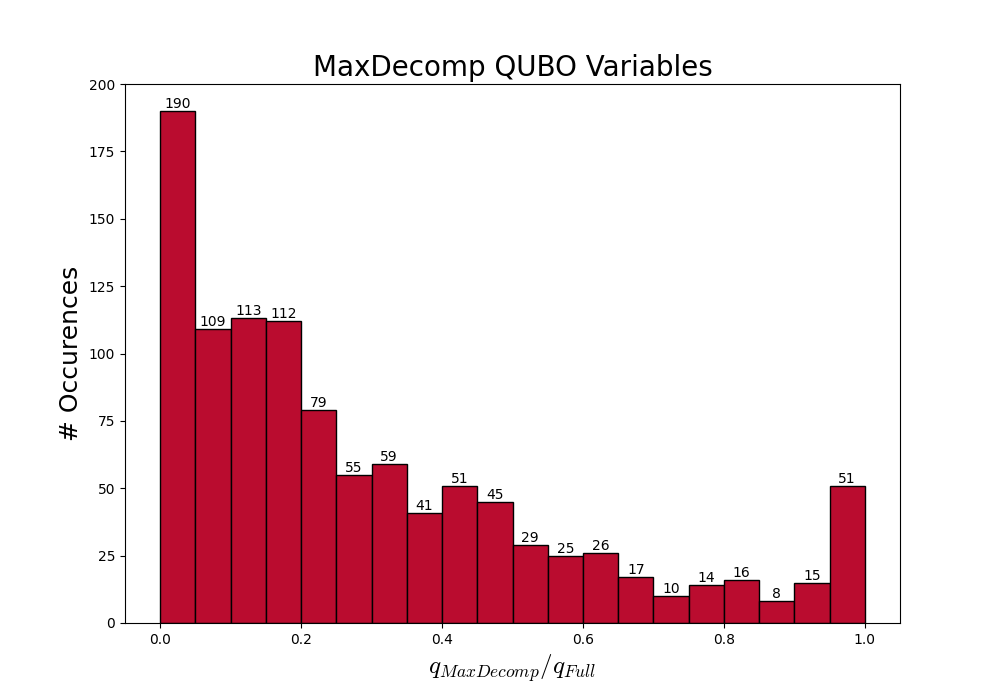}
    \label{fig:part-dist}
    \caption{Distributions of $q_{Reduced}$ (left) and $q_{MaxDecomp}$ (right) as a proportion of $q_{Full}$ within the MIPLIB 2017 Collection Set.} 
\end{figure}

To study how much the Reduced and Decomposed forms of our formulations decrease the number of variables in our QUBOs in practice, we calculate values of $\nu$, $\mu$, and the size of the largest decomposition, which we refer to as MaxDecomp, for each problem in the MIPLIB 2017 collection. We then use these values to determine how many variables each formulation would require as a percent of the Full formulation. 
On average, the Reduced form requires 32\% of the number of QUBO variables, while the MaxDecomp form requires 30\% of the number of variables compared to the Full formulation. The full distribution of the percent of $\pi,\sigma$ entries needed is shown in Figure \ref{fig:part-dist} .

When $\frac{q_{Reduced}}{q_{Full}}$ is closer to 1, more of the MIP variables and constraints could be reasonably permuted. However, the overall size of the QUBO will remain quite large. On the other hand, if $\frac{q_{Reduced}}{q_{Full}}$ is closer to 0, then the QUBO will be much smaller, but the MIP likely has very little symmetry worth exploiting. Therefore, when working with an MIP problem, the trade-off of QUBO size and potential for symmetry must be balanced. 

\begin{figure}
    \centering
    \includegraphics[width=0.49\linewidth]{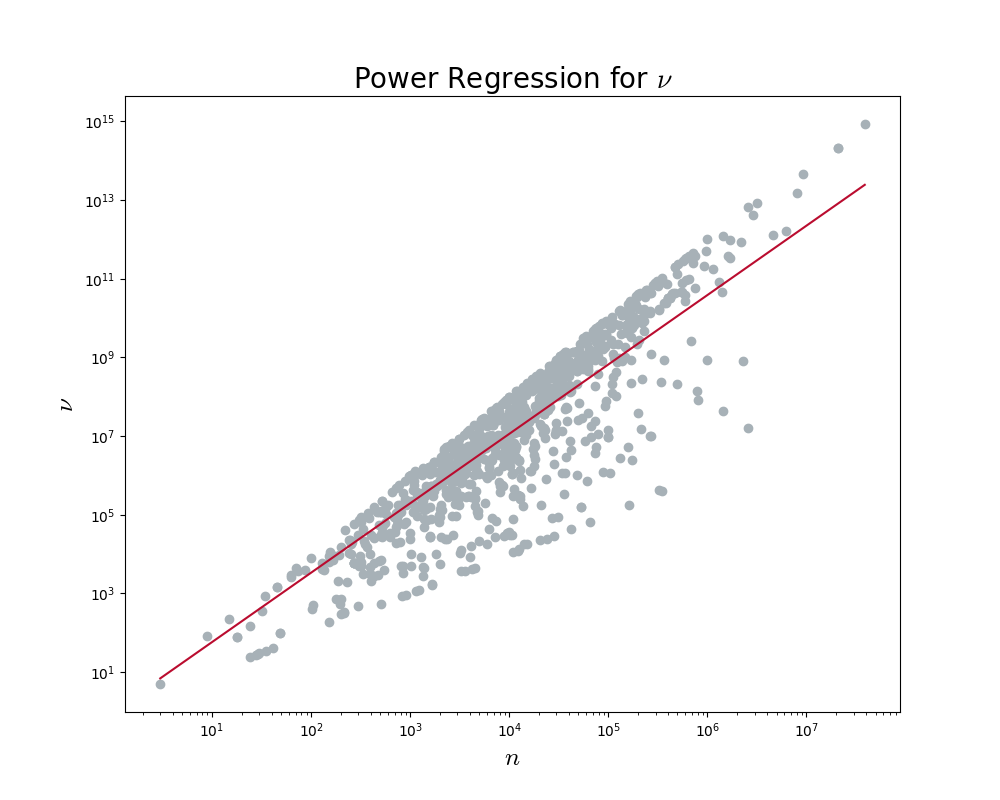}
    \includegraphics[width=0.49\linewidth]{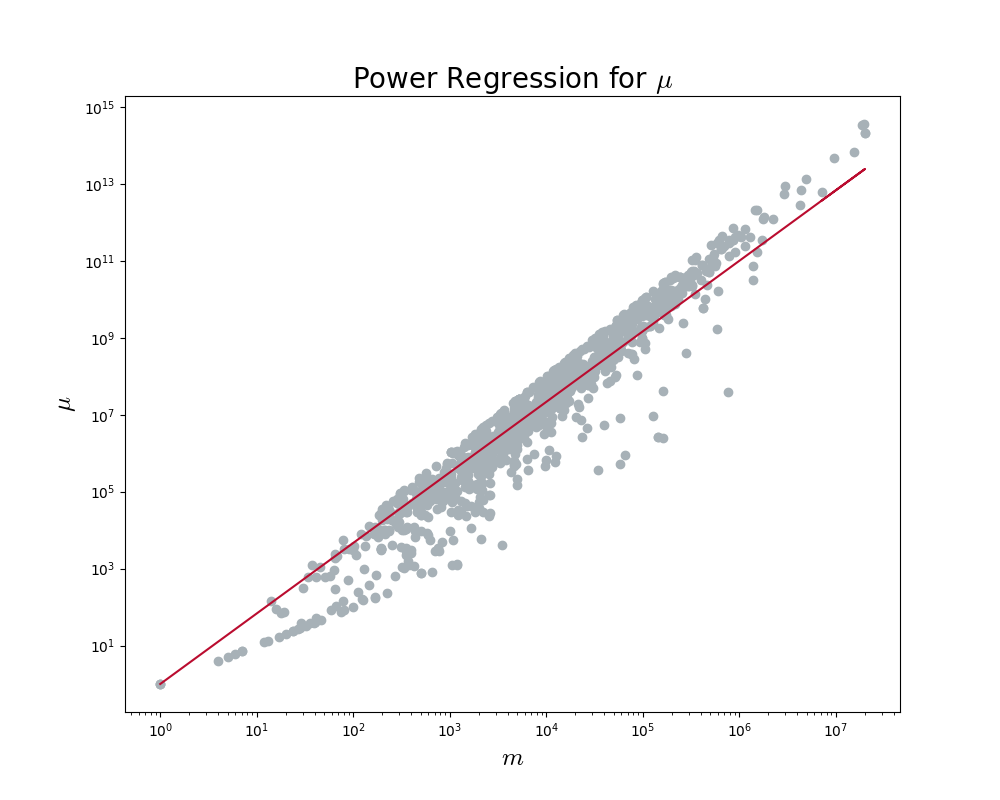}
    \caption{Power regressions of $\nu$ and $\mu$ as functions of $n$ and $m$ within the MIPLIB 2017 Collection Set.}
    \label{fig:regressions}
\end{figure}

We have also calculated a power regression of the form $y=x^k$ to estimate $\nu$ and $\mu$ as a function of $n$ and $m$ in the test set. We have $\nu \approx n^{1.764}$ and $\mu \approx m^{1.834}$. Scatter-plots of the data can be seen in Figure \ref{fig:regressions}.

\subsection{DWave Embeddings}

To use a DWave quantum annealer, the QUBO problem graph must be embedded on the hardware's topology. While we can use Proposition \ref{prop:qubits} to determine the number of qubits needed in a Zephyr graph (with $t=4$), this is a rather loose bound in practice. Thus we apply the \texttt{find\_embedding} routine from the \texttt{minorminer} library provided with Dwave software, which is a heuristic for embedding a source graph (our QUBO) to a target graph (an appropriately sized Zephyr graph). For each reduced QUBO, we set the target graph as $Z_g$ where $g=\frac{\mu+\mu+8}{16}$, the graph for which a $K_{\nu+\mu}$ graph could be embedded. 

Due to the computation time of the heuristic routine, we were only able to test the embedding heuristic on the 12 smallest embeddings. For the largest four instances within this set, no embedding was found within the routine's time limit. Complete results can be seen in Table \ref{tab:dwave-embeddings}. On average, the heuristic required 39.8\% the amount of qubits that the complete graph embedding requires. 

\begin{figure}
    \centering
    \includegraphics[width=0.49\linewidth]{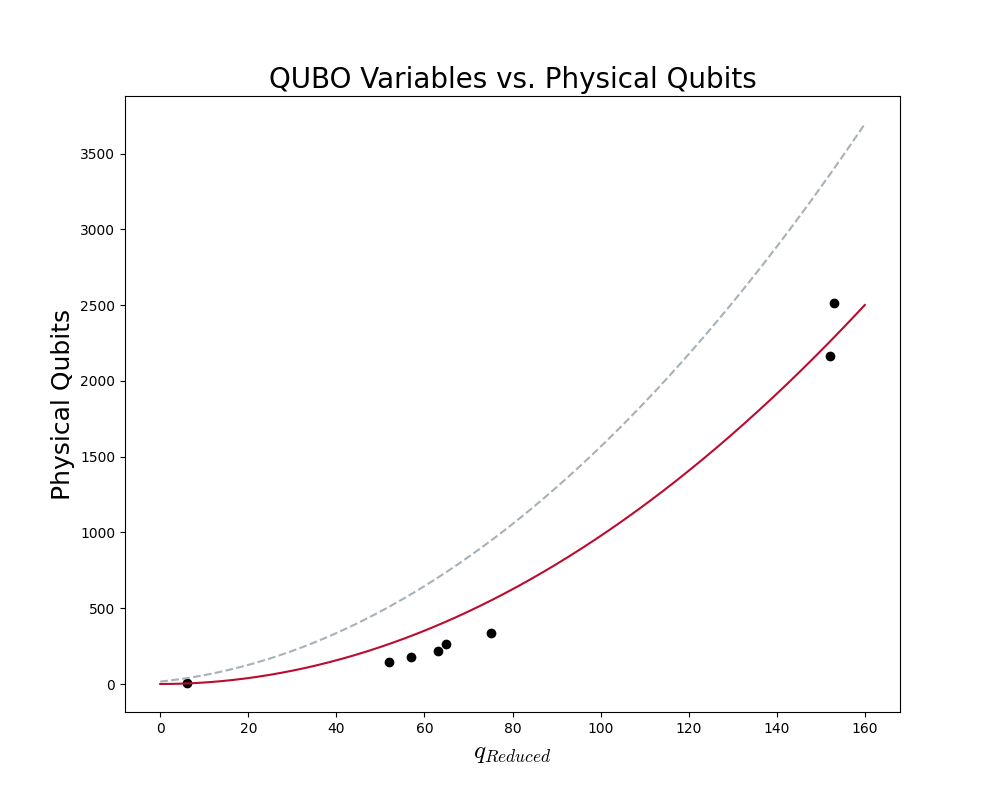}
    \includegraphics[width=0.49\linewidth]{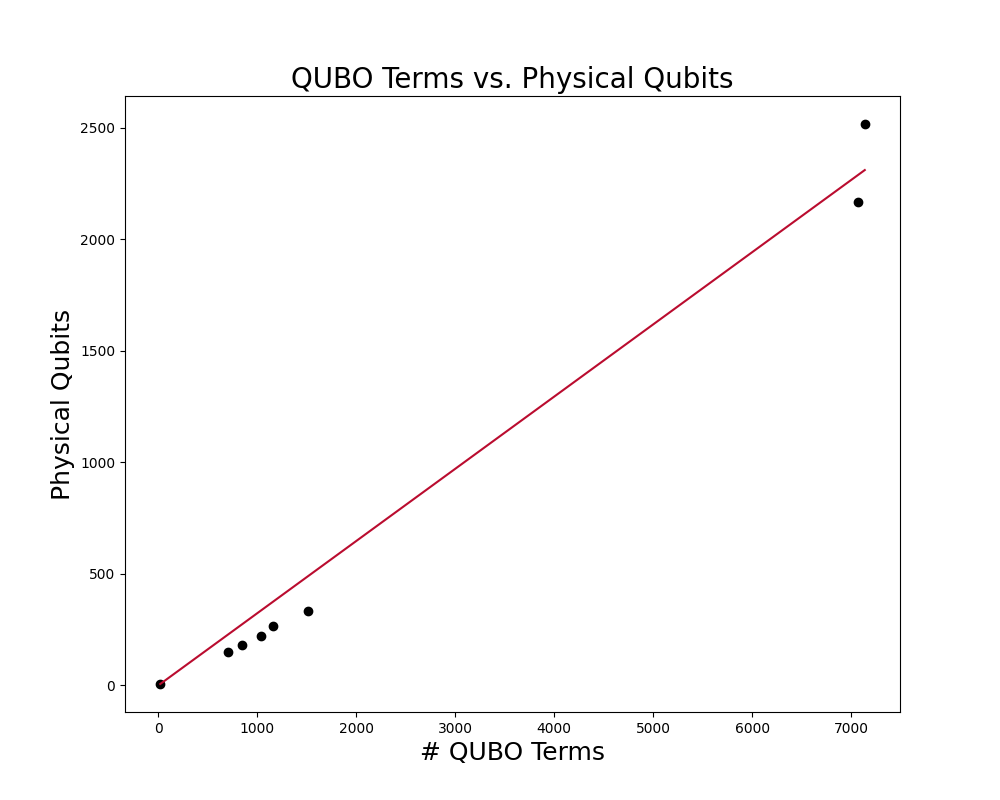}
    \caption{Regressions of physical qubits required for embedding on Zephyr graphs vs $q_{Reduced}$ (left), where the dashed line is the bound described in Proposition \ref{prop:qubits}, and the number of QUBO terms (right)}
    \label{fig:qubit-regressions}
\end{figure}

We run a regression on the number of physical qubits required to embed each problem on a Zephyr graph, finding that $qubits \approx 0.98q_{Reduced}^2$.  A stronger correlation was found with the number of QUBO terms, which is reasonable as the connectivity (or lack thereof) is what necessitates the embeddings in the first place. We found that $qubits \approx 0.324(\#terms)$. Visualizations of these regressions can be seen in Figure \ref{fig:qubit-regressions}. Based on the formula for the number of qubits required to embed the $K_{q_{Reduced}}$ graph as well as the fact that there are $O(v^2)$ edges on a graph with $v$ vertices, it is not surprising that the required physical qubits are still $O(q_{Reduced}^2)$ and $O(\#terms)$. Combining our regressions, for the reduced symmmetry formulation of an MIP with $n$ variables and $m$ constraints, we would expect to require a Zephyr topoplogy quantum annealer with $0.98(n^{1.764}+m^{1.834})^2 \in O(n^{3.528}+m^{3.668})$ qubits. 

\pagebreak

\bibliographystyle{splncs04}
\bibliography{references}

@article{abbas2024challenges,
  title={Challenges and opportunities in quantum optimization},
  author={Abbas, Amira and Ambainis, Andris and Augustino, Brandon and B{\"a}rtschi, Andreas and Buhrman, Harry and Coffrin, Carleton and Cortiana, Giorgio and Dunjko, Vedran and Egger, Daniel J and Elmegreen, Bruce G and others},
  journal={Nature Reviews Physics},
  pages={1--18},
  year={2024},
  publisher={Nature Publishing Group}
}

@article{lubinski2024optimization,
  title={Optimization applications as quantum performance benchmarks},
  author={Lubinski, Thomas and Coffrin, Carleton and McGeoch, Catherine and Sathe, Pratik and Apanavicius, Joshua and Bernal Neira, David and Quantum Economic Development Consortium and others},
  journal={ACM Transactions on Quantum Computing},
  volume={5},
  number={3},
  pages={1--44},
  year={2024},
  publisher={ACM New York, NY}
}

@article{proctor2025benchmarking,
  title={Benchmarking quantum computers},
  author={Proctor, Timothy and Young, Kevin and Baczewski, Andrew D and Blume-Kohout, Robin},
  journal={Nature Reviews Physics},
  volume={7},
  number={2},
  pages={105--118},
  year={2025},
  publisher={Nature Publishing Group UK London}
}

@article{nannicini2019performance,
  title={Performance of hybrid quantum-classical variational heuristics for combinatorial optimization},
  author={Nannicini, Giacomo},
  journal={Physical Review E},
  volume={99},
  number={1},
  pages={013304},
  year={2019},
  publisher={APS}
}

@inproceedings{aaronson2016space,
  title={The space" just above" BQP},
  author={Aaronson, Scott and Bouland, Adam and Fitzsimons, Joseph and Lee, Mitchell},
  booktitle={Proceedings of the 2016 ACM Conference on Innovations in Theoretical Computer Science},
  pages={271--280},
  year={2016}
}

@article{jun2023hubo,
  title={HUBO and QUBO models for prime factorization},
  author={Jun, Kyungtaek and Lee, Hyunju},
  journal={Scientific Reports},
  volume={13},
  number={1},
  pages={10080},
  year={2023},
  publisher={Nature Publishing Group UK London}
}

@article{ajagekar2022hybrid,
  title={Hybrid classical-quantum optimization techniques for solving mixed-integer programming problems in production scheduling},
  author={Ajagekar, Akshay and Al Hamoud, Kumail and You, Fengqi},
  journal={IEEE Transactions on Quantum Engineering},
  volume={3},
  pages={1--16},
  year={2022},
  publisher={IEEE}
}

@inproceedings{babai2018group,
  title={Group, graphs, algorithms: the graph isomorphism problem},
  author={Babai, L{\'a}szl{\'o}},
  booktitle={Proceedings of the International Congress of Mathematicians: Rio de Janeiro 2018},
  pages={3319--3336},
  year={2018},
  organization={World Scientific}
}

@article{mohammadisiahroudi2025improvements,
  title={Improvements to quantum interior point method for linear optimization},
  author={Mohammadisiahroudi, Mohammadhossein and Wu, Zeguan and Augustino, Brandon and Carr, Arielle and Terlaky, Tam{\'a}s},
  journal={ACM Transactions on Quantum Computing},
  volume={6},
  number={1},
  pages={1--24},
  year={2025},
  publisher={ACM New York, NY}
}

@article{nannicini2024fast,
  title={Fast quantum subroutines for the simplex method},
  author={Nannicini, Giacomo},
  journal={Operations Research},
  volume={72},
  number={2},
  pages={763--780},
  year={2024},
  publisher={INFORMS}
}

@article{montanaro2020quantum,
  title={Quantum speedup of branch-and-bound algorithms},
  author={Montanaro, Ashley},
  journal={Physical Review Research},
  volume={2},
  number={1},
  pages={013056},
  year={2020},
  publisher={APS}
}

@misc{woerner2021solving,
  title={Solving mixed integer optimization problems on a hybrid classical-quantum computing system},
  author={Woerner, Stefan and Nannicini, Giacomo and Barkoutsos, Panagiotis and Tavernelli, Ivano},
  year={2021},
  month=jan # "~26",
  publisher={Google Patents},
  note={US Patent 10,902,085}
}

@article{brown2024copositive,
  title={A copositive framework for analysis of hybrid ising-classical algorithms},
  author={Brown, Robin and Bernal Neira, David E and Venturelli, Davide and Pavone, Marco},
  journal={SIAM Journal on Optimization},
  volume={34},
  number={2},
  pages={1455--1489},
  year={2024},
  publisher={SIAM}
}

@inproceedings{wei2024hybrid,
  title={Hybrid quantum-classical computing via Dantzig-Wolfe decomposition for integer linear programming},
  author={Wei, Xinliang and Liu, Jiyao and Fan, Lei and Guo, Yuanxiong and Han, Zhu and Wang, Yu},
  booktitle={2024 33rd International Conference on Computer Communications and Networks (ICCCN)},
  pages={1--9},
  year={2024},
  organization={IEEE}
}

@article{ellinas2024hybrid,
  title={A hybrid quantum--classical algorithm for mixed-integer optimization in power systems},
  author={Ellinas, Petros and Chevalier, Samuel and Chatzivasileiadis, Spyros},
  journal={Electric Power Systems Research},
  volume={235},
  pages={110835},
  year={2024},
  publisher={Elsevier}
}

@article{paterakis2023hybrid,
  title={Hybrid quantum-classical multi-cut benders approach with a power system application},
  author={Paterakis, Nikolaos G},
  journal={Computers \& Chemical Engineering},
  volume={172},
  pages={108161},
  year={2023},
  publisher={Elsevier}
}

@inproceedings{zhao2022hybrid,
  title={Hybrid quantum benders’ decomposition for mixed-integer linear programming},
  author={Zhao, Zhongqi and Fan, Lei and Han, Zhu},
  booktitle={2022 IEEE Wireless Communications and Networking Conference (WCNC)},
  pages={2536--2540},
  year={2022},
  organization={IEEE}
}

@article{svensson2023hybrid,
  title={Hybrid quantum-classical heuristic to solve large-scale integer linear programs},
  author={Svensson, Marika and Andersson, Martin and Gr{\"o}nkvist, Mattias and Vikst{\aa}l, Pontus and Dubhashi, Devdatt and Ferrini, Giulia and Johansson, G{\"o}ran},
  journal={Physical Review Applied},
  volume={20},
  number={3},
  pages={034062},
  year={2023},
  publisher={APS}
}

@article{junger2021quantum,
  title={Quantum annealing versus digital computing: An experimental comparison},
  author={J{\"u}nger, Michael and Lobe, Elisabeth and Mutzel, Petra and Reinelt, Gerhard and Rendl, Franz and Rinaldi, Giovanni and Stollenwerk, Tobias},
  journal={Journal of Experimental Algorithmics (JEA)},
  volume={26},
  pages={1--30},
  year={2021},
  publisher={ACM New York, NY, USA}
}

@article{carlson2012miso,
  title={MISO unlocks billions in savings through the application of operations research for energy and ancillary services markets},
  author={Carlson, Brian and Chen, Yonghong and Hong, Mingguo and Jones, Roy and Larson, Kevin and Ma, Xingwang and Nieuwesteeg, Peter and Song, Haili and Sperry, Kimberly and Tackett, Matthew and others},
  journal={Interfaces},
  volume={42},
  number={1},
  pages={58--73},
  year={2012},
  publisher={INFORMS}
}

@inproceedings{boudot2020comparing,
  title={Comparing the difficulty of factorization and discrete logarithm: a 240-digit experiment},
  author={Boudot, Fabrice and Gaudry, Pierrick and Guillevic, Aurore and Heninger, Nadia and Thom{\'e}, Emmanuel and Zimmermann, Paul},
  booktitle={Annual International Cryptology Conference},
  pages={62--91},
  year={2020},
  organization={Springer}
}

@article{peng2019factoring,
  title={Factoring larger integers with fewer qubits via quantum annealing with optimized parameters},
  author={Peng, WangChun and Wang, BaoNan and Hu, Feng and Wang, YunJiang and Fang, XianJin and Chen, XingYuan and Wang, Chao},
  journal={SCIENCE CHINA Physics, Mechanics \& Astronomy},
  volume={62},
  number={6},
  pages={60311},
  year={2019},
  publisher={Springer}
}

@article{willsch2023large,
  title={Large-scale simulation of Shor’s quantum factoring algorithm},
  author={Willsch, Dennis and Willsch, Madita and Jin, Fengping and De Raedt, Hans and Michielsen, Kristel},
  journal={Mathematics},
  volume={11},
  number={19},
  pages={4222},
  year={2023},
  publisher={Multidisciplinary Digital Publishing Institute}
}

@article{quinton2025quantum,
  title={Quantum annealing applications, challenges and limitations for optimisation problems compared to classical solvers},
  author={Quinton, Finley Alexander and Myhr, Per Arne Sevle and Barani, Mostafa and Crespo del Granado, Pedro and Zhang, Hongyu},
  journal={Scientific Reports},
  volume={15},
  number={1},
  pages={12733},
  year={2025},
  publisher={Nature Publishing Group UK London}
}

@article{jain2011qip,
  title={Qip= pspace},
  author={Jain, Rahul and Ji, Zhengfeng and Upadhyay, Sarvagya and Watrous, John},
  journal={Journal of the ACM (JACM)},
  volume={58},
  number={6},
  pages={1--27},
  year={2011},
  publisher={ACM New York, NY, USA}
}

@article{margot2009symmetry,
  title={Symmetry in integer linear programming},
  author={Margot, Fran{\c{c}}ois},
  journal={50 Years of Integer Programming 1958-2008: From the Early Years to the State-of-the-Art},
  pages={647--686},
  year={2009},
  publisher={Springer}
}

@article{ostrowski2015modified,
  title={Modified orbital branching for structured symmetry with an application to unit commitment},
  author={Ostrowski, James and Anjos, Miguel F and Vannelli, Anthony},
  journal={Mathematical Programming},
  volume={150},
  pages={99--129},
  year={2015},
  publisher={Springer}
}

@article{pfetsch2019computational,
  title={A computational comparison of symmetry handling methods for mixed integer programs},
  author={Pfetsch, Marc E and Rehn, Thomas},
  journal={Mathematical Programming Computation},
  volume={11},
  pages={37--93},
  year={2019},
  publisher={Springer}
}

@article{bodi2013algorithms,
  title={Algorithms for highly symmetric linear and integer programs},
  author={B{\"o}di, Richard and Herr, Katrin and Joswig, Michael},
  journal={Mathematical Programming},
  volume={137},
  number={1-2},
  pages={65--90},
  year={2013},
  publisher={Springer}
}

@article{mckay2007nauty,
  title={Nauty user’s guide (version 2.4)},
  author={McKay, Brendan D},
  journal={Computer Science Dept., Australian National University},
  pages={225--239},
  year={2007}
}

@inproceedings{aaronson2010bqp,
  title={BQP and the polynomial hierarchy},
  author={Aaronson, Scott},
  booktitle={Proceedings of the forty-second ACM symposium on Theory of computing},
  pages={141--150},
  year={2010}
}

@article{hua2020improved,
  title={Improved QUBO formulation of the graph isomorphism problem},
  author={Hua, Richard and Dinneen, Michael J},
  journal={SN Computer Science},
  volume={1},
  pages={1--18},
  year={2020},
  publisher={Springer}
}

@article{calude2017qubo,
  title={QUBO formulations for the graph isomorphism problem and related problems},
  author={Calude, Cristian S and Dinneen, Michael J and Hua, Richard},
  journal={Theoretical Computer Science},
  volume={701},
  pages={54--69},
  year={2017},
  publisher={Elsevier}
}

@article{Lucas:2014,
  author          = {Lucas, Andrew},
  journal         = {Frontiers in Physics},
  number          = {5},
  title           = {Ising formulations of many NP problems},
  volume          = {2},
  year            = {2014},
  DOI             = {https://10.3389/fphy.2014.00005}
}

@inproceedings{junttila2007engineering,
  title={Engineering an efficient canonical labeling tool for large and sparse graphs},
  author={Junttila, Tommi and Kaski, Petteri},
  booktitle={2007 Proceedings of the Ninth Workshop on Algorithm Engineering and Experiments (ALENEX)},
  pages={135--149},
  year={2007},
  organization={SIAM}
}

@INPROCEEDINGS{shor1994shors,

  author={Shor, P.W.},

  booktitle={Proceedings 35th Annual Symposium on Foundations of Computer Science}, 

  title={Algorithms for quantum computation: discrete logarithms and factoring}, 

  year={1994},

  volume={},

  number={},

  pages={124-134},

  doi={10.1109/SFCS.1994.365700}}

@article{glover2022quantum1,
  title={Quantum bridge analytics I: a tutorial on formulating and using QUBO models},
  author={Glover, Fred and Kochenberger, Gary and Hennig, Rick and Du, Yu},
  journal={Annals of Operations Research},
  volume={314},
  number={1},
  pages={141--183},
  year={2022},
  publisher={Springer}
}

@article{liberti2012reformulations,
  title={Reformulations in mathematical programming: automatic symmetry detection and exploitation},
  author={Liberti, Leo},
  journal={Mathematical Programming},
  volume={131},
  pages={273--304},
  year={2012},
  publisher={Springer}
}

@article{Abbott2019anneal,
   title={A hybrid quantum-classical paradigm to mitigate embedding costs in quantum annealing},
   volume={17},
   ISSN={1793-6918},
   url={http://dx.doi.org/10.1142/S0219749919500424},
   DOI={10.1142/s0219749919500424},
   number={05},
   journal={International Journal of Quantum Information},
   publisher={World Scientific Pub Co Pte Ltd},
   author={Abbott, Alastair A. and Calude, Cristian S. and Dinneen, Michael J. and Hua, Richard},
   year={2019},
   month=aug, pages={1950042} }

@article{glover2022quantum2,
  title={Quantum Bridge Analytics II: QUBO-Plus, network optimization and combinatorial chaining for asset exchange},
  author={Glover, Fred and Kochenberger, Gary and Ma, Moses and Du, Yu},
  journal={Annals of Operations Research},
  volume={314},
  number={1},
  pages={185--212},
  year={2022},
  publisher={Springer}
}

@article{fiedler1988doubly,
  title={Doubly stochastic matrices and optimization},
  author={Fielder, M},
  journal={Mathematical research},
  volume={45},
  pages={44--51},
  year={1988}
}

@techreport{boothby2021zephyr,
    author = {Boothby, Kelly and King, Andrew D. and Raymond, Jack},
    title = {Zephyr Topology of D-Wave Quantum Processors},
    institution = {DWave},
    year = {2021}
}

@article{boothby2016fast,
  title={Fast clique minor generation in Chimera qubit connectivity graphs},
  author={Boothby, Tomas and King, Andrew D and Roy, Aidan},
  journal={Quantum Information Processing},
  volume={15},
  pages={495--508},
  year={2016},
  publisher={Springer}
}

@article{bolusani2024scip,
  title={The SCIP optimization suite 9.0},
  author={Bolusani, Suresh and Besan{\c{c}}on, Mathieu and Bestuzheva, Ksenia and Chmiela, Antonia and Dion{\'\i}sio, Jo{\~a}o and Donkiewicz, Tim and van Doornmalen, Jasper and Eifler, Leon and Ghannam, Mohammed and Gleixner, Ambros and others},
  journal={arXiv preprint arXiv:2402.17702},
  year={2024}
}

@article{wang2025rbm,
  title={RBM-Based Simulated Quantum Annealing for Graph Isomorphism Problems},
  author={Wang, Yukun and Shen, Yingtong and Zhang, Zhichao and Wan, Linchun},
  journal={arXiv preprint arXiv:2503.07749},
  year={2025}
}

@article{developers2020d,
  title={D-Wave Hybrid Solver Service: An Overview},
  author={Developers, D-Wave},
  journal={D-Wave Systems Inc., Tech. Rep},
  year={2020}
}

\pagebreak
\appendix
\section{Results of Heuristic Embeddings}
\begin{table}[]
\centering
\caption{Results of embedding our Reduced Symmtery Detecting QUBO on the smallest MIPLIB 2017 instances on Zephyr topology with $t=4$}
\resizebox{\textwidth}{!}{%
\begin{tabular}{@{}lllllllllll@{}}
\toprule
Problem   & $n$ & $m$ & $\nu$ & $\mu$ & $q_{Reduced}$ & \begin{tabular}[c]{@{}l@{}}\# QUBO \\ Terms\end{tabular} & $g$ & \begin{tabular}[c]{@{}l@{}}Qubits\\ (Heuristic)\end{tabular} & \begin{tabular}[c]{@{}l@{}}Qubits\\ $K_{\nu+\mu}$\end{tabular} & Proportion \\ \midrule
ej        & 3   & 1   & 5     & 1     & 6             & 18                                                       & 1   & 7                                                            & 39                                                             & 18\%       \\
flugpl    & 18  & 18  & 80    & 72    & 152           & 7066                                                     & 10  & 2166                                                         & 3360                                                           & 64\%       \\
flugplinf & 18  & 19  & 80    & 73    & 153           & 7139                                                     & 11  & 2517                                                         & 3402                                                           & 74\%       \\
gen-ip016 & 28  & 24  & 28    & 24    & 52            & 706                                                      & 4   & 148                                                          & 510                                                            & 29\%       \\
gen-ip036 & 29  & 46  & 29    & 46    & 75            & 1516                                                     & 6   & 334                                                          & 945                                                            & 35\%       \\
gen-ip054 & 30  & 27  & 30    & 27    & 57            & 843                                                      & 5   & 180                                                          & 594                                                            & 30\%       \\
gen-ip021 & 35  & 28  & 35    & 28    & 63            & 1036                                                     & 5   & 220                                                          & 702                                                            & 31\%       \\
gen-ip002 & 41  & 24  & 41    & 24    & 65            & 1161                                                     & 5   & 265                                                          & 740                                                            & 36\%       \\ \bottomrule
\end{tabular}%
}
\label{tab:dwave-embeddings}
\end{table}
\end{document}